\documentclass[12pt,reqno]{amsart}
\usepackage{pkg}

\newcommand{\AKMaxcut}{MR1488983}
\newcommand{\DemboMontanariIsing}{MR2650042}
\newcommand{\DemboMontanariGibbs}{MR2643563}
\newcommand{\DMS}{arXiv:1110.4821}
\newcommand{\DyerFriezeJerrumSIAM}{MR1936657}
\newcommand{\DyerFriezeJerrumFOCS}{MR1917561}
\newcommand{\DommersEtAl}{MR2733399}
\newcommand{\GalanisEtAl}{springerlink:10.1007/978364222935048}
\newcommand{\JansonEtAl}{MR1782847}
\newcommand{\JerrumSinclair}{MR1237164}
\newcommand{\LiLuYin}{arXiv:1111.7064}
\newcommand{\LubyVigodaRSA}{MR1716763}

\newcommand{\MMS}{springerlink:10.1007/s00440-010-0315-6}
\newcommand{\MWW}{MR2475668}
\newcommand{\SinclairSrivastavaThurley}{arXiv:1107.2368}
\newcommand{\SlyFOCS}{10.1109/FOCS.2010.34}
\newcommand{\Weitz}{MR2277139}
\newcommand{\Zachary}{MR714953}

\newcommand{\lwc}{\to_{\mathrm{\it loc}}}
\newcommand{\lpc}{\to_{\mathrm{\it loc}}}
\newcommand{\rt}{o}
\newcommand{\tree}{T}
\newcommand{\trees}{\mathcal{T}}
\newcommand{\treespr}{{\trees_{\edge}}}
\newcommand{\treereg}{\T}

\newcommand{\biptreereg}{\mathrm{\textbf{\textup T}}}
\newcommand{\bipP}{\mathrm{\textbf{\textup P}}}
\newcommand{\bipE}{\mathrm{\textbf{\textup E}}}
\newcommand{\col}{\tau}

\newcommand{\match}{\mathfrak{m}}

\newcommand{\fer}{\mathrm{\textup f}}
\newcommand{\af}{\mathrm{\textup af}}
\newcommand{\free}{\circ}
\newcommand{\BP}{{\mathrm{\textsf{BP}}}}
\newcommand{\BPF}{{\mathrm{\textsf{F}}}}
\newcommand{\BPG}{{\mathrm{\textsf{G}}}}
\newcommand{\danglb}[1]{\llbracket #1 \rrbracket}
\newcommand{\gibbs}{\mathscr{G}}
\newcommand{\simplex}{\De}
\newcommand{\spins}{\mathscr{X}}
\newcommand{\vpsi}{{\bar\psi}}
\newcommand{\vxi}{{\bar\xi}}
\newcommand{\vertex}{{\mathrm{vx}}}
\newcommand{\edge}{{\mathrm{e}}}

\newcommand{\NP}{{\mathrm{\textsc{np}}}}
\newcommand{\RP}{{\mathrm{\textsc{rp}}}}

\newcommand{\pras}{{\mathrm{\textsc{pras}}}}
\newcommand{\fptas}{{\mathrm{\textsc{fptas}}}}
\newcommand{\fpras}{{\mathrm{\textsc{fpras}}}}
\newcommand{\maxcut}{{\mathrm{\textsc{max-cut}}}}
\newcommand{\cut}{\mathrm{cut}}
\newcommand{\OUT}{\mathrm{out}}
\newcommand{\IN}{\mathrm{in}}

\begin{document}

\title[Hardness of counting for two-spin models on $d$-regular graphs]{The computational hardness of counting \\ in two-spin models on $d$-regular graphs}

\author[A.\ Sly]{$^*$Allan Sly}
\address{Department of Statistics, University of California, Berkeley
\newline\indent Evans Hall, Berkeley, California 94720}

\author[N.\ Sun]{$^\dagger$Nike Sun}
\address{Department of Statistics, Stanford University
\newline\indent Sequoia Hall, 390 Serra Mall, Stanford, California 94305}

\date{\today}

\thanks{$^*$Research partially supported by Alfred P.\ Sloan Research Fellowship. \newline\indent $^\dagger$Research partially supported by Department of Defense NDSEG Fellowship.}

\maketitle

\vspace{-1cm}

\begin{abstract}
The class of two-spin systems contains several important models, including random independent sets and the Ising model of statistical physics. We show that for both the hard-core (independent set) model and the anti-ferromagnetic Ising model with arbitrary external field, it is $\NP$-hard to approximate the partition function or approximately sample from the model on $d$-regular graphs when the model has non-uniqueness on the $d$-regular tree. Together with results of Jerrum--Sinclair, Weitz, and Sinclair--Srivastava--Thurley giving $\fpras$'s for all other two-spin systems except at the uniqueness threshold, this gives an almost complete classification of the computational complexity of two-spin systems on bounded-degree graphs.

Our proof establishes that the normalized log-partition function of any two-spin system on bipartite locally tree-like graphs converges to a limiting ``free energy density'' which coincides with the (non-rigorous) Bethe prediction of statistical physics. We use this result to characterize the local structure of two-spin systems on locally tree-like bipartite expander graphs, which then become the basic gadgets in a randomized reduction to approximate $\maxcut$. Our approach is novel in that it makes no use of the second moment method employed in previous works on these questions.
\end{abstract}

\section{Introduction}

Spin systems are stochastic models defined by local interactions on networks. The class of spin systems includes well-known combinatorial counting and constraint satisfaction problems. In this paper we classify the complexity of approximating the partition function for all homogeneous two-spin systems on bounded-degree graphs.

When interactions favor agreement of adjacent spins, the model is said to be \emph{ferromagnetic}. Jerrum and Sinclair~\cite{\JerrumSinclair} gave a fully polynomial-time randomized approximation scheme ($\fpras$) for approximating the partition function (the normalizing constant in the probability distribution) of the ferromagnetic Ising model, which covers all ferromagnetic two-spin systems. For anti-ferromagnetic systems such as the hard-core and anti-ferromagnetic Ising models, the complexity of approximating the partition function depends on the model parameters, and is known to be $\NP$-hard when the interactions are sufficiently strong. Our first main result establishes that the computational transition for such models on $d$-regular graphs is located precisely at the uniqueness threshold (see Defn.~\ref{d:unique}) for the corresponding model on the $d$-regular tree.

\bThm\label{t:hc.comp}
For $d\ge3$ and $\lm>\lm_c(d)= \frac{(d-1)^{d-1}}{(d-2)^d}$, unless $\NP=\RP$ there exists no $\fpras$ for the partition function of the hard-core model with fugacity $\lm$ on $d$-regular graphs.
\eThm

The transition point $\lm_c(d)$ is the uniqueness threshold for the hard-core model on the $d$-regular tree: it marks the point above which distant boundary conditions have a non-vanishing influence on the spin at the root. In a seminal paper~\cite{\Weitz}, Weitz used computational tree methods to provide a $\fptas$ for the partition function of the hard-core model on graphs of maximum degree $d$ at any $\lm<\lm_c(d)$. Together with Weitz's result, Thm.~\ref{t:hc.comp} completes the classification of the complexity of the hard-core model except at the threshold~$\lm_c$.

Previously it was shown that there is no $\fpras$ for the hard-core model at $\lambda d\ge10000$~\cite{\LubyVigodaRSA}. In the case of $\lm=1$ this was improved to $d\ge25$~\cite{\DyerFriezeJerrumFOCS,\DyerFriezeJerrumSIAM}, using random regular bipartite graphs as basic gadgets in a hardness reduction. Mossel et al.~\cite{\MWW} showed that local {\sc mcmc} algorithms are exponentially slow for $\lm>\lm_c(d)$, and conjectured that $\lm_c$ is in fact the threshold for existence of an $\fpras$.

The first rigorous result establishing a computational transition at the uniqueness threshold appeared in~\cite{\SlyFOCS}, where hardness was shown for $\lm_c(d)<\lm<\lm_c(d)+\ep(d)$ for some $\epsilon(d)>0$. The proof relies on a detailed analysis of the hard-core model on random bipartite graphs, which are then used in a randomized reduction to $\maxcut$. More precisely the result of~\cite{\SlyFOCS} gives hardness subject to a technical condition which was an artifact of a difficult second moment calculation from~\cite{\MWW}, and which could only be verified for $\lm<\lm_c(d)+\ep(d)$. Hardness was subsequently shown by Galanis et al.~\cite{\GalanisEtAl} for all $\lm>\lm_c(d)$ when $d\ne4,5$ by verifying the technical condition of~\cite{\SlyFOCS}.

In this paper we follow a different approach which is more conceptual and completely circumvents second moment method calculations. Moreover the same method of proof gives the analogous result for anti-ferromagnetic Ising models with arbitrary external field:

\bThm\label{t:ising.comp}
For $d\ge3$, $B\in\R$ and $\be<\be_{c,\af}(B,d)<0$, unless $\NP=\RP$ there does not exist an $\fpras$ for the partition function of the anti-ferromagnetic Ising model with inverse temperature $\be$ and external field $B$ on $d$-regular graphs.
\eThm

Here $\be_{c,\af}(B,d)$ denotes the uniqueness threshold for the anti-ferromagnetic Ising model with external field $B$ on the $d$-regular tree. Extending the methods of Weitz~\cite{\Weitz}, Sinclair et al.~\cite{\SinclairSrivastavaThurley} (see also \cite{\LiLuYin}) gave a $\fptas$ for the anti-ferromagnetic Ising model on $d$-regular graphs at inverse temperature $\be>\beta_{c,\af}(B,d)$, so together with Thm.~\ref{t:ising.comp} this again establishes that the computational transition coincides with the tree uniqueness threshold.

The hard-core and anti-ferromagnetic Ising models together encompass all (non-degenerate) homogeneous two-spin systems on $d$-regular graphs (see \S\ref{ss:reduc.two}). Thus, the results of~\cite{\Weitz,\JerrumSinclair,\SinclairSrivastavaThurley} combined with Thms.~\ref{t:hc.comp} and~\ref{t:ising.comp} give a full classification of the computational complexity of approximating the partition function for (homogeneous) two-spin systems on $d$-regular graphs, except at the uniqueness thresholds $\lambda_c(d)$ and $\beta_{c,\af}(B,d)$.

In fact, we will show inapproximability in non-uniqueness regimes in a strong sense: not only does there not exist an $\fpras$, but for any fixed choice of model parameters and $d$ there exists $c>0$ such that it is $\NP$-hard even to approximate the partition function within a factor of $e^{cn}$ on the class of $d$-regular graphs.

\subsection*{Independent results of Galanis--\v{S}tefankovi\v{c}--Vigoda}
In a simultaneous and independent work, Galanis, \v{S}tefankovi\v{c} and Vigoda~\cite{GalanisStefankovicVigoda} established the result of Thm.~\ref{t:hc.comp}, and Thm.~\ref{t:ising.comp} in the case of zero external field ($B=0$). Their methods differ from ours: they analyze the second moment of the partition function on random bipartite $d$-regular graphs, and establish the condition necessary to apply the approach of~\cite{\SlyFOCS}. Their proof analyzes a difficult optimization of a real function in several variables by relating the problem to certain tree recursions.

\subsection{Reduction to $\maxcut$ via bipartite graphs}

Our proof is based on a detailed characterization (Thm.~\ref{t:loc}) of the local structure of anti-ferromagnetic two-spin systems on symmetric bipartite $d$-regular locally tree-like graphs. Specifically, we show that the joint distribution of all the spins in a large neighborhood of a typical vertex in the graph converges to a known (Gibbs) measure on the $d$-regular tree. Under the additional assumption that that the graph is an edge expander, when the model has non-uniqueness on the $d$-regular tree the spin distribution on the graph is divided into $+$ and $-$ phases where one or the other side of the graph has a linear number more vertices with $+$ spin.

Our main results Thms.~\ref{t:hc.comp}~and~\ref{t:ising.comp} are then proved by a variation on the construction of~\cite{\SlyFOCS}, using the bipartite graphs in a randomized reduction approximate $\maxcut$ on $3$-regular graphs, which is known to be $\NP$-hard~\cite{\AKMaxcut}. First, we use Thm.~\ref{t:loc} to construct a symmetric bipartite $d$-regular locally tree-like graph $G$ of large constant size such that, conditioned on the phase of the global configuration, spins at distant vertices are asymptotically independent with known marginals depending only on the side of the graph (Propn.~\ref{p:product}).

Given a $3$-regular graph $H$ on which we wish to approximate $\maxcut$, first we take a disjoint copy $G_v$ of $G$ for each vertex $v\in H$. After removing $3k$ edges from each $G_v$, for each edge $(u,v)\in H$ we add $k$ edges joining each side of $G_u$ to the corresponding side of $G_v$ in such a way that the resulting graph $H^G$ is $d$-regular.

The connections between gadgets do not substantially change the spin distributions inside them, and in particular the $\pm$ phases remain. The anti-ferromagnetic nature of the interation, however, results in neighboring copies of $G$ in $H^G$ preferring to be in opposing phases. Using the asymptotic conditional independence result Propn.~\ref{p:product} we can estimate the partition function for the model on $H^G$ restricted to configurations of given phase on each copy of $G$ within a factor of $e^{\ep\abs{H}}$ (Lem.~\ref{l:HGcutProb}). We find that the distribution is concentrated on configurations where the vector of phases gives a good cut of $H$, and the effect is strengthened as $k$ is increased. Thus, for any $\ep>0$, by taking $k$ (hence $G$) to be sufficiently large a $(1+\ep)$-approximation of $\maxcut(H)$ can be determined from the partition function of the model on $H^G$, thereby completing the reduction.

Our reduction depends crucially on the detailed picture of the spin distribution developed in Thm.~\ref{t:loc} and Propn.~\ref{p:product}. Using methods developed in~\cite{\MMS}, these results in turn are obtained as consequences of precise asymptotics for the partition function of two-spin models on bipartite $d$-regular graphs: we show that the log-partition function, normalized by the number of vertices in the graph, has an asymptotic value, the ``free energy density,'' which is easily computed from the (non-rigorous) ``Bethe prediction'' of statistical physics (see \S\ref{ss:bethe}). This is a result of independent interest, since lower bounds for partition functions on graphs have proved to be in general challenging. Asymptotics for the partition function on general tree-like graphs were established for the ferromagnetic Ising model in~\cite{\DemboMontanariIsing,\DommersEtAl,\DMS}, and for more general spin systems in uniqueness regimes in~\cite{\DMS}. Our result for anti-ferromagnetic models is stated somewhat informally as follows; for the precise statement see Thm.~\ref{t:z}.

\bThm\label{t:intro.z}
For any non-degenerate homogeneous two-spin model on bipartite $d$-regular locally tree-like graphs, the log-partition function normalized by the number of vertices has an asymptotic value which coincides with the Bethe free energy prediction.
\eThm

In the remainder of this introductory section we formally introduce the models which we consider. We then define the notion of local (weak) convergence of graphs and give precise statements of our results on the partition function (Thm.~\ref{t:z}) and local structure (Thm.~\ref{t:loc}) of these models on bipartite graphs.

\subsection{Definition of spin systems}

Let $G=(V,E)$ be a finite undirected graph, and $\spins$ a finite alphabet of \emph{spins}. A \emph{spin system} or \emph{spin model} on $G$ is a probability measure on the space of \emph{(spin) configurations} $\vec\si\in\spins^V$ of form
\beq\label{e:fm}
\nu^{\vec\psi}_G(\vec\si)
= \f{1}{Z_G(\vec\psi)}
	\prod_{(ij)\in E} \psi(\si_i,\si_j)
	\prod_{i\in V} \vpsi(\si_i),
\eeq
where $\psi$ is a symmetric function $\spins^2\to\R_{\ge0}$, $\vpsi$ is a positive function $\spins\to\R_{\ge0}$, and $Z_G(\vec\psi)$ is the normalizing constant, called the \emph{partition function}. The pair $\vec\psi\equiv(\psi,\vpsi)$ is called a \emph{specification} for the spin system~\eqref{e:fm}.

In this paper we consider spin systems with an alphabet of size two; without loss $\spins\equiv\set{\pm1}$. The \emph{Ising} model on $G$ at inverse temperature $\be$ and external field $B$ is given by
\beq\label{e:ising}
\nu^{\be,B}_G(\vec\si)
=\f{1}{Z_G(\be,B)}
	\prod_{(ij)\in E} e^{\be\si_i\si_j}
	\prod_{i\in V} e^{B\si_i}.
\eeq
The \emph{hard-core} (or \emph{independent set}) model on $G$ at activity or fugacity $\lm$ is given by
\beq\label{e:hardcore}
\nu^\lm_G(\vec\si)
=\f{1}{Z_G(\lm)}
	\prod_{(ij)\in E} \Ind{\bar\si_i\bar\si_j\ne1}
	\prod_{i\in V} \lm^{\bar\si_i}
\eeq
where $\bar\si\equiv\Ind{\si=+1}=(1+\si)/2$. The edge interaction has no temperature parameter and includes a hard constraint. Our definition \eqref{e:hardcore} is trivially equivalent to the standard definition of the hard-core model which has spin $0$ in place of $-1$, but we take $\spins=\set{\pm1}$ throughout to unify the notation.

\subsection{Local convergence and the Bethe prediction}
\label{ss:intro.z}

If $G$ is any graph and $v$ a vertex in $G$, write $\ball{t}{v}$ for the subgraph induced by the vertices of $G$ at graph distance at most $t$ from $v$, and $\pd v\equiv\ball{1}{v}\setminus\set{v}$ for the neighbors of $v$. We let $\tree\equiv(\tree,\rt)$ denote a general tree with root $\rt$, with $\tree^t\equiv\ball{t}{\rt}\subseteq\tree$ the subtree of depth $t$. We also fix $d$ throughout and write $\treereg\equiv(\treereg,\rt)$ for the rooted $d$-regular tree.

\bdfn\label{d:lpc}
Let $G_n=(V_n=[n],E_n)$ be a sequence of (random) finite undirected graphs, and let $I_n$ denote a uniformly random vertex in $V_n$. The sequence $G_n$ is said to \emph{converge locally} to the $d$-regular tree $\treereg$ if for all $t\ge0$, $\ball{t}{I_n}$ converges to $\treereg^t$ in distribution with respect to the joint law $\P_n$ of $(G_n,I_n)$: that is, $\lim_{n\to\infty}\P_n(\ball{t}{I_n}\cong\treereg^t)=1$ (where $\cong$ denotes graph isomorphism).
\edfn

We write $\E_n$ for expectation with respect to $\P_n$ and impose the following integrability condition on the degree of $I_n$:

\bdfn
The sequence $G_n$ is \emph{uniformly sparse} if the random variables $\abs{\pd I_n}$ are uniformly integrable, that is, if
$$\lim_{L\to\infty}\limsup_{n\to\infty}\E_n[ \abs{\pd I_n} \Ind{ \abs{\pd I_n}\ge L }]=0.$$
\edfn

We assume throughout that $G_n$ ($n\ge1$) is a uniformly sparse graph sequence converging locally to the $d$-regular tree $\treereg$; this setting is hereafter denoted $G_n\lpc\treereg$. The \emph{free energy density} for a specification $\vec\psi$ on $G_n$ is defined by
\beq\label{e:lim}
\phi
\equiv\lim_{n\to\infty}\phi_n
\equiv\lim_{n\to\infty} \f{1}{n}\E_n[\log Z_n],\quad
Z_n\equiv Z_{G_n}(\vec\psi),
\eeq
provided the limit exists. For ferromagnetic spin systems on a broad class of locally tree-like graphs, heuristic methods from statistical physics yield an explicit (conjectural) formula for the value of $\phi$, the so-called ``Bethe prediction'' $\Phi$ whose definition we recall in \S\ref{ss:bethe}. For anti-ferromagnetic two-spin models, the Bethe prediction is well-defined only on graph sequences $G_n$ which are nearly bipartite, in the following sense: let $\treereg_+$ denote the $d$-regular tree $\treereg$ with vertices colored $+1$ (black) or $-1$ (white) according to whether they are at even or odd distance from the root $\rt$; let $\treereg_-$ be $\treereg_+$ with the colors reversed. Let $\biptreereg$ be the random tree which equals $\treereg_+$ or $\treereg_-$ with equal probability; write $\bipP$ for the law of $\biptreereg$ and $\bipE$ for expectation with respect to $\bipP$.

\bdfn\label{d:lwc}
For $G_n\lpc\treereg$, we say the $G_n$ are \emph{nearly bipartite}, and write $G_n\lwc\biptreereg$ (equivalently $G_n\lwc\bipP$), if there exists a (not necessarily proper) black-white coloring of $G_n$ such that for all $t\ge0$, $\ball{t}{I_n}\to\biptreereg^t$ in distribution.
\edfn

The precise statement of Thm.~\ref{t:intro.z} is then as follows:

\bThm\label{t:z}
Let $\vec\psi$ specify a non-degenerate homogeneous two-spin system.
\bnm[(a)]
\item If $\vec\psi$ is ferromagnetic, then $\phi$ exists for any $G_n\lpc\treereg$ and equals $\Phi_{\set{\treereg}}$ as defined by \eqref{e:Bethe} (and given more explicitly by \eqref{e:bethe.f}).
\item If $\vec\psi$ is anti-ferromagnetic, then $\phi$ exists for any $G_n\lwc\biptreereg$ and equals $\Phi_{\set{\treereg_\pm}}$ as defined in \eqref{e:Bethe} (and given more explicitly by \eqref{e:bethe.af}).
\enm
\eThm

\brmk\label{r:trees}
Hereafter we treat $G_n\lpc\treereg$ and $G_n\lwc\biptreereg$ in a unified manner when possible by writing $G_n\lwc\P_\trees$ for $\P_\trees$ the uniform measure on $\trees$, which always denotes either $\set{\treereg}$ or $\set{\treereg_\pm}$. We write $\E_\trees$ for expectation with respect to $\P_\trees$.
\ermk

\subsection{Local structure of measures}
\label{ss:intro.loc}

Under some additional assumptions on $G_n$, Thm.~\ref{t:z}, together with the arguments of~\cite{\MMS}, characterizes the asymptotic local structure of the spin systems $\nu_n\equiv\nu_{G_n}$. For $G_n\lwc\biptreereg$, let $\col:V_n\to\set{\pm}$ denote the given black-white coloring of the vertices of $G_n$ (hereafter writing $\pm$ as shorthand for $\pm1$). We say that $G_n$ is \emph{symmetric} if it is isomorphism-invariant to reversing the black-white coloring. For a spin configuration $\vec\si\in G_n$ we define the \emph{phase} of $\vec\si$ to be
$$Y(\vec\si)\equiv\sgn \sum_i\col_i\si_i,\quad
\text{where}\quad\sgn x\equiv\Ind{x\ge0}-\Ind{x<0}.
$$
Let $\nu^\pm_n$ denote the measure $\nu_n$ conditioned on the configurations of $\pm$ phase: that is,
$$\nu^\pm_n(\vec\si)
\equiv \f{1}{Z^\pm_n}
	\Ind{Y(\vec\si)=\pm}
	\prod_{(ij)\in E_n}\psi(\si_i,\si_j)
	\prod_{i\in V_n}\vpsi(\si_i),$$
where $Z^\pm_n$ is the partition function restricted to the $\pm$ configurations. We will characterize the local structure of the measures $\nu^\pm_n$ on graph sequences satisfying an edge-expansion assumption, as follows:

\bdfn\label{d:exp}
A graph $G=(V,E)$ is a \emph{$(\de,\gam,\lm)$-edge expander} if, for any set of vertices $S\subseteq V$ with $\de\abs{V}\le\abs{S}\le\gam\abs{V}$, there are at least $\lm\abs{S}$ edges joining $S$ to $V\setminus S$.
\edfn

The measures $\nu^\pm_n$ will be related to Gibbs measures on the infinite tree. In particular, recall the definition of (Gibbs) uniqueness:

\bdfn\label{d:unique}
For a rooted tree $\tree$, let $\gibbs_\tree$ denote the set of Gibbs measures for the specification $\vec\psi$ on $\tree$. The specification is said to have \emph{(Gibbs) uniqueness} (on $\tree$) if $\abs{\gibbs_\tree}=1$.
\edfn

Recalling Rmk.~\ref{r:trees}, let $\gibbs_\trees$ denote the space of mappings $\nu:\tree\mapsto\nu(\tree)$, $\tree\in\trees$ (with $\gibbs_{\set{\treereg}}\hookrightarrow\gibbs_{\set{\treereg_\pm}}$ in the obvious manner). When $\trees=\set{\treereg_\pm}$ we write $\nu_\pm$ as shorthand for $\nu(\treereg_\pm)$.

\bdfn\label{d:trans}
An element $\nu\in\gibbs_\trees$ is \emph{translation-invariant} if for $(\tree,\rt)\in\trees$ and any vertex $x\in \tree$, the law on spin configurations of $(\tree,x)$ induced by $\nu(\tree,\rt)$ coincides with $\nu(\tree,x)$.\footnote{If $\trees=\set{\treereg}$ this agrees with the usual definition of translation-invariance, whereas if $\trees=\set{\treereg_\pm}$ then the projections $\nu(\treereg^\pm)$ are \emph{semi}-translation-invariant.}
\edfn

For a two-spin model, let  $\nu^+$ (resp.\ $\nu^-$) be the elements of $\gibbs_\trees$ defined by conditioning on all spins identically equal to $1$ on the $t$-th level of black (resp.\ white) vertices and taking the weak limit as $t\to\infty$; the $\nu^\pm$ are translation-invariant. The projections $\mu^+\equiv\nu^+_+\equiv\nu^+(\treereg_+)$ and $\mu^-\equiv\nu^-_+\equiv\nu^-(\treereg_+)$, disregarding the black-white coloring on $\treereg_+$, are the extremal semi-translation-invariant Gibbs measures for the model on $\treereg$, and by symmetry
$$\mu^+=\nu^-_-\equiv\nu^-(\treereg_-),\quad
\mu^-=\nu^+_-\equiv\nu^+(\treereg_-).$$
The model has uniqueness if and only if $\mu^+=\mu^-$.

\bdfn\label{d:lwc.msr}
For $G_n\sim\P_n$ a random graph sequence and $\nu_n$ any law on spin configurations $\vec\si_n$ of $G_n$, we say that $\P_n\otimes\nu_n$ \emph{converges locally (weakly)} to $\P_\trees\otimes\nu$ (for $\nu\in\gibbs_\trees$), and write $\P_n\otimes\nu_n\lwc\P_\trees\otimes\nu$, if it holds for all $t\ge0$ that $(\ball{t}{I_n},\vec\si_{\ball{t}{I_n}})$ converges in distribution to $(\tree^t,\vec\si_t)$ where $\tree\sim\P_\trees$ and $\vec\si_t$ is the restriction to $\tree^t$ of $\vec\si\sim\nu(\tree)$.
\edfn

\brmk
In \cite[Defn.~2.3]{\MMS} three forms $A,B,C$ of local convergence of measures are distinguished, with $C\Rightarrow B\Rightarrow A$. Our Defn.~\ref{d:lwc.msr} corresponds to the weakest form $A$: however, as explained in the proof of \cite[Thm.~2.4~(II)]{\MMS}, if the $(\nu(\tree))_{\tree\in\trees}$ are \emph{extremal} Gibbs measures then $A,B,C$ are easily seen to be equivalent, so convergence in the sense of Defn.~\ref{d:lwc.msr} implies convergence in the {\it a priori} stronger sense of
$$\tvd{\P_n[(\ball{t}{I_n},\vec\si_{\ball{t}{I_n}})=\cdot] - \P_\trees[(\tree^t,\vec\si_{\tree^t})=\cdot]}\to0.$$
\ermk

\bThm\label{t:loc}
For any anti-ferromagnetic two-spin system on $G_n\lwc\biptreereg$, the following hold:
\bnm[(a)]
\item If the $G_n$ are symmetric, then $\P_n\otimes\nu_n\lwc\bipP\otimes[(\nu^+ + \nu^-)/2]$.
\item If for all $\de>0$ the $G_n$ are $(\de,1/2,\lm_\de)$-edge expanders for some $\lm_\de>0$, then
	\beq\label{e:locConvergence}
    \P_n\otimes\nu^\pm_n\lwc\bipP\otimes\nu^\pm.
    \eeq
    Further, with $\angl{\,}_\mu$ denoting expectation with respect to the Gibbs measure $\mu$,
    \beq\label{e:magnetization}
        \frac1n Y(\vec\si) \sum_{i \in V} \col_i\si_i
        \to
        \frac12[
        \angl{\si_\rt}_{\mu^+}
        -\angl{\si_\rt}_{\mu^-}
        ]\quad\text{in probability.}
    \eeq
\enm
\eThm

\subsection*{Outline of the paper}

In \S\ref{s:z} we review the Bethe prediction in the $d$-regular setting and prove Thm.~\ref{t:intro.z} (in its form Thm.~\ref{t:z}). In \S\ref{s:loc} we show how to deduce Thm.~\ref{t:loc} from Thm.~\ref{t:z} by the methods of~\cite{\MMS}. In \S\ref{s:comp} we prove the approximate conditional independence statement (Propn.~\ref{p:product}) and demonstrate the randomized reduction to $\maxcut$ to prove our main results Thms.~\ref{t:hc.comp}~and~\ref{t:ising.comp}.

\section{Partition function for two-spin models}
\label{s:z}

In this section we prove Thm.~\ref{t:z}, establishing the free energy density $\phi$ (and verifying the Bethe prediction) for two-spin models on graph sequences $G_n\lwc\biptreereg$. We refer to~\cite{\DemboMontanariGibbs,\DMS} for more general background and references on the Bethe prediction, and in \S\ref{ss:bethe} describe only its specialization to the $d$-regular setting. In \S\ref{ss:reduc.two} we show that for purposes of computing $\phi$ on $d$-regular locally tree-like graph sequences, all non-degenerate two-spin systems reduce to Ising or hard-core. In \S\ref{ss:interpolate} we compute the free energy density for these models by applying an interpolation scheme described in~\cite{\DMS}, thereby completing the proof of Thm.~\ref{t:z}.

\subsection{The Bethe prediction}
\label{ss:bethe}

Recalling the notation of Rmk.~\ref{r:trees}, we now review the Bethe prediction for $G_n\lwc\P_\trees$. Given $\trees$, let $\treespr$ denote the set of trees $\tree$ rooted not at a vertex but at an oriented edge $x\to y$, obtained by distinguishing an oriented edge in $\tree\in\trees$ and forgetting the root. Elements of $\trees,\treespr$ are regarded modulo isomorphism: thus if
$\trees=\set{\treereg}$ then $\treespr=\set{(\treereg,\rt\to j)}$, and if $\trees=\set{\treereg_\pm}$ then $\treespr=\set{(\treereg_\pm,\rt\to j)}$.

Let $\simplex$ denote the $(\abs{\spins}-1)$-dimensional simplex of probability measures on $\spins$. A \emph{message} is a mapping $h:\treespr\to\simplex$; we write $\cH\equiv\cH(\trees)$ for the space of messages on $\treespr$. For $T\in\trees$, $x\to y$ in $T$, and $h\in\cH$, write $h_{x\to y}$ for the image of $(T,x\to y)\in\treespr$ under $h$, and define
$$\Phi_T(h)\equiv\Phi^\vertex_T(h)-\Phi^\edge_T(h)$$
where
\begin{align*}
\Phi^\vertex_T(h)&\equiv
	\log\bigg\{ \sum_{\si_\rt}\vpsi(\si_\rt)
		\prod_{j\in\pd\rt}\bigg(
			\sum_{\si_j}\psi(\si_\rt,\si_j)
				h_{j\to\rt}(\si_j)
			\bigg) \bigg\},\\
\Phi^\edge_T(h)&\equiv \f12
	\sum_{j\in\pd\rt} \log\bigg\{
		\sum_{\si_\rt,\si_j} \psi(\si_\rt,\si_j)
			h_{\rt\to j}(\si_\rt)
			h_{j\to\rt}(\si_j)
		\bigg\}.
\end{align*}
The \emph{Bethe free energy functional} on $\cH(\trees)$ is defined by $\Phi_\trees(h)\equiv\E_\trees[\Phi_T(h)]$.

The \emph{Bethe} or \emph{belief propagation (BP) recursion} is the map
$$\BP\equiv\BP_\trees:\cH(\trees)\to\cH(\trees),\quad
(\BP h)_{x\to y}(\si)
	\equiv\bar\BPF[(h_{v\to x})_{v\in\pd x\setminus y}]$$
for $\bar\BPF:\simplex^{d-1}\to\simplex$ defined by
\beq\label{e:bpf}
[\bar\BPF(\vec h)](\si)
\cong \vpsi(\si) \prod_{j=1}^{d-1} \bigg\{
		\sum_{\si_j}\psi(\si,\si_j) h_j(\si_j)
		\bigg\},\quad \vec h\equiv(h_1,\ldots,h_{d-1}) \in\simplex^{d-1}
\eeq
(where $\cong$ denotes equivalence up to a positive normalizing factor).

\bdfn
For any homogeneous spin system on $G_n\lpc\P_\trees$, the \emph{Bethe prediction} is that the free energy density $\phi$ of \eqref{e:lim} exists and equals
\beq\label{e:Bethe}
\Phi\equiv\Phi_\trees\equiv\sup_{h\in\cH_\star}\Phi_\trees(h)
\eeq
with $\cH_\star\equiv\cH_\star(\trees)\subseteq\cH(\trees)$ the set of all fixed points of $\BP_\trees$.
\edfn

For $h\in\simplex$ write $\BPF(h)\equiv\bar\BPF(h,\ldots,h)$: then $\cH_\star(\set{\treereg})$ corresponds simply to the fixed points of $\BPF$ in simplex. For $h\in\cH(\set{\treereg_\pm})$ we write $h_\pm\equiv h(\treereg_\pm,\rt\to j)\in\simplex$: then any $h\in\cH_\star(\set{\treereg_\pm})$ must satisfy $h_\pm=\BPF(h_\mp)$, so $\cH_\star(\set{\treereg_\pm})$ corresponds to the fixed points of the \emph{double} recursion $\BPF^{(2)}\equiv\BPF\circ\BPF$.

In verifying the Bethe prediction we will identify the fixed points attaining the supremum in \eqref{e:Bethe}. In the anti-ferromagnetic case, with $h^+$ (resp.\ $h^-$) denoting the elements $h\in\cH_\star(\set{\treereg_\pm})$ maximizing $h_+(+)$ (resp.\ $h_-(+)$), we will see that
\beq\label{e:bethe.af}
\Phi_{\set{\treereg_\pm}} = \Phi_{\set{\treereg_\pm}}(h^+) = \Phi_{\set{\treereg_\pm}}(h^-).
\eeq
Explicitly, $h^+_+=h^-_-$ (resp.\ $h^+_-=h^-_+$) will be the fixed points of $\BPF^{(2)}$ giving maximal (resp.\ minimal) probability to spin $+$. The ferromagnetic case reduces to the Ising model: here, with $h^\pm$ denoting the elements of $\cH_\star(\set{\treereg})$ maximizing $h_{\rt\to j}(\pm)$ on $\treereg$, we will see that
\beq\label{e:bethe.f}
\Phi_{\set{\treereg}}=\Phi_{\set{\treereg}}(h^{\sgn B}).
\eeq
The remainder of this section is devoted to the proof of Thm.~\ref{t:z}.

\subsection{Reduction to Ising and hard-core on $d$-regular graphs}
\label{ss:reduc.two}

We first show that for the computation of the free energy density, all (non-degenerate) homogeneous two-spin models on graph sequences $G_n\lwc\biptreereg$ reduce to either the Ising or hard-core model. Indeed, let $\vec\psi\equiv(\psi,\vpsi)$ be a specification for a two-spin system with alphabet $\spins=\set{\pm}$. If we define $\vec\psi'$ by $\psi'(\si,\si')\equiv\psi(\si,\si')\vpsi(\si)^{1/d}\vpsi(\si')^{1/d},$ and $\vpsi'(\si)\equiv1$, then
$$\f{1}{n}\log Z_G(\vec\psi)
-\f{1}{n}\log Z_G(\vec\psi')
= O(\E_n[ \abs{\pd I_n} \Ind{ \abs{\pd I_n}\ne d}]),$$
which for $G_n\lpc\treereg$ tends to zero as $n\to\infty$ by uniform sparsity. Therefore we assume without loss $\vpsi\equiv1$, and consider the possibilities for $\psi$:
\bnm
\item[(1)] If $\psi>0$, then
$\psi(\si,\si')= e^{B_0} e^{\be\si\si'} e^{B\si/d} e^{B\si'/d}$ for $\be,B,B_0$ defined by
$$\f{\psi(+,+)}{\psi(-,-)}=e^{4B/d},\quad
\f{\psi(+,+)\psi(-,-)}{\psi(+,-)^2}=e^{4\be},\quad
\psi(+,+)\psi(+,-)^2\psi(-,-)=e^{4B_0},$$
so $\phi_n-(d/2)B_0$ is asymptotically equal to the free energy density for the Ising model on $G_n$ with parameters $(\be,B)$.

\item[(2)] If $\psi(+,-)=\psi(-,+)>0$ and $\psi(-,-)>\psi(+,+)=0$, then, recalling $\bar\si\equiv\Ind{\si=+}$, we have $\psi(\si,\si')=e^{B_0} \Ind{\bar\si\bar\si'\ne1} \lm^{\bar\si/d} \lm^{\bar\si'/d}$ for $B_0,\lm$ defined by
$$\psi(-,-)\equiv e^{B_0},\quad
\f{\psi(+,-)}{\psi(-,-)}\equiv \lm^{1/d}.$$
Therefore $\phi_n-(d/2)B_0$ is asymptotically equal to the free energy density for the independent set model on $G_n$ at fugacity $\lm$.
\enm
The remaining two-spin models are degenerate, with free energy density which is easy to calculate:
\bnm
\item[(3)] Suppose $\psi(+,-)=\psi(-,+)=0$, so that $\psi(\si,\si')$ may be written as $\Ind{\si=\si'} e^{B_0} e^{B\si/d}e^{B\si'/d}$. Then
$$\phi_n = B_0 \f{\E_n[\abs{E_n}]}{n}
	+ B + \f{1}{n}\E_n\bigg[\sum_{j=1}^{k(G_n)} \log (1 + e^{-2 B\abs{C_j}})\bigg]$$
where the sum is taken over the connected components $C_1,\ldots,C_{k(G_n)}$ of $G_n$. We claim $\phi_n\to\phi=(d/2)B_0+B$: we have $\liminf_{n\to\infty} (\phi_n-\phi)\ge0$ (using uniform sparsity), and
$$\limsup_{n\to\infty}\,(\phi_n-\phi)
\le\limsup_{n\to\infty}\log 2
	\f{\E_n[k(G_n)]}{n},
$$
so it suffices to show $\E_n[k(G_n)]/n\to0$. Indeed, if this fails then there exists $\ep>0$ such that for infinitely many $n$, the event $\set{k(G_n)\ge\ep n}$ occurs with $\P_n$-probability at least $\ep$. On this event, $G_n$ has at least $\ep n/2$ components of size $\le 2/\ep$, so for $t>\log_k(2/\ep)$, $\limsup_{n\to\infty}\P_n(\ball{t}{I_n}\not\cong\treereg^t)\ge \ep^2/2>0$, in contradiction of $G_n\lpc\treereg$.

\item[(4)] Suppose instead $\psi(+,+)=\psi(-,-)=0$ while $\psi(+,-)=\psi(-,+)>0$. If the $G_n$ are not exactly bipartite then $\phi_n=-\infty$. If they are exactly bipartite then
$$\phi_n = \log\psi(+,-) \f{\E_n[\abs{E_n}]}{n}
	+ \log 2\f{\E_n[k(G_n)]}{n},$$
and by the observation of (3) this converges to $\phi=(d/2)\log\psi(+,-)$.

\enm

\subsection{Bethe interpolation}
\label{ss:interpolate}

We now evaluate the hard-core and Ising free energy densities by interpolating in the model parameters. Write $\xi\equiv\log\psi$, $\vxi\equiv\log\vpsi$, and for the hard-core model take $B\equiv\log\lm$. Let $\angl{\,}^{\be,B}_n$ denote expectation with respect to $\nu^{\be,B}_n\equiv\nu^{\be,B}_{G_n}$, and define
\begin{align*}
a^\vertex_n(\be,B)
&\equiv\pd_B\phi_n(\be,B)
	= \E_n[\angl{ \pd_B\vxi(\si_{I_n}) }^{\be,B}_n],\\
a^\edge_n(\be,B)
&\equiv\pd_\be\phi_n(\be,B)
	= \f12\E_n\Big[
		\sum_{j\in\pd I_n}
			\angl{\pd_\be\xi(\si_{I_n},\si_j)}^{\be,B}_n
		\Big]
\end{align*}
(with $a^\edge_n(\be,B)\equiv0$ for hard-core). We also define analogous quantities on the limiting tree $\tree\sim\P_\trees$: for $h\in\cH$ let
\begin{align}
\nonumber
a^\vertex(\be,B,h)
&\equiv a^\vertex_\trees(\be,B,h)
\equiv \E_\trees[ \danglb{\pd_B\vxi(\si_\rt)}^{\be,B}_h ],\\
a^\edge(\be,B,h)
&\equiv a^\edge_\trees(\be,B,h)
\equiv \f12 \E_\trees\Big[
		\sum_{j\in\pd\rt} \danglb{ \pd_\be\xi(\si_\rt,\si_j)}^{\be,B}_h
		\Big]
\label{e:a.tree}
\end{align}
where $\danglb{\cdot}_h$ denotes expectation with respect to the measure $\nu^h_{\tree^1}$ on spin configurations on $\tree^1$ defined by
$$\nu^h_{\tree^1}(\vec\si_{\tree^1}=\cdot)
\cong
\vpsi(\si_\rt) \prod_{j\in\pd\rt} \psi(\si_\rt,\si_j) h_{j\to\rt}(\si_j).$$
The following lemma, describing our interpolation scheme, may be verified directly or obtained as a consequence of \cite[Propn.~2.4]{\DMS}. We always interpolate in one parameter at a time, keeping the other fixed and suppressing it from the notation.

\blem\label{l:int}
If for $B\in[B_0,B_1]$ we have $h\equiv h(B)\in\cH^B_\star$ which is continuous and of bounded total variation in $B$, then
$$\Phi_\trees(B_1)-\Phi_\trees(B_0)
= \int_{B_0}^{B_1} a^\vertex_\trees(B,h) \d B.$$
The same result holds for $B,a^\vertex_\trees$ replaced with $\be,a^\edge_\trees$.
\elem

We now make explicit the connection between BP fixed points and Gibbs measures; for a discussion in a more general setting and further references see \cite[Rmk.~2.6]{\DMS}. Recall that for $\tree\in\trees$, $\gibbs_\tree$ denotes the set of Gibbs measures for the specification $\vec\psi$ on $\tree$, and $\gibbs_\trees$ denotes the space of mappings $\tree\mapsto\nu(\tree)\in\gibbs_\tree$ with $\tree\in\trees$. An element $\mu=\nu(\tree)\in\gibbs_\tree$ is a \emph{Markov chain} or \emph{splitting Gibbs measure} (see~\cite{\Zachary}) if there exists a collection $h^\mu\equiv(h^\mu_{x\to y})$ of elements of $\simplex$ indexed by the oriented edges of $\tree$ such that for any finite connected induced subgraph $U=(V_U,E_U)$ of $\tree$,
\beq\label{e:markov}
\mu(\vec\si_U)
= \f{1}{z} \prod_{i\in V_U}\vpsi(\si_i)
	\prod_{(ij)\in E_U}\psi(\si_i,\si_j)
	\prod_{j\in\pd U}
	\bigg\{
	\sum_{\si_j}
	\psi(\si_{p(j)},\si_j) h^\mu_{j\to p(j)}(\si_j)
	\bigg\},
\eeq
where $p(j)$ denotes the unique neighbor of $j$ inside $U$ for $j$ belonging to the external boundary $\pd U$ of $U$. Extremal Gibbs measures are Markov chains but the converse is false. We say that an element $\nu\in\gibbs_\trees$ is \emph{Markovian} if $\nu(\tree)$ is a Markov chain for each $\tree\in\trees$: the associated collection $h^\nu\equiv(h^{\nu(\tree)})_{\tree\in\trees}$ is called an \emph{entrance law}: it satisfies consistency conditions imposed by \eqref{e:markov} (which closely resemble the BP equation), and the correspondence between Markovian $\nu\in\gibbs_\trees$ and entrance laws $h^\nu$ is bijective. If $\nu$ is also \emph{translation-invariant} (in the sense of Defn.~\ref{d:trans}), then each $h^{\nu(\tree)}_{x\to y}$ depends only on the isomorphism class of $(\tree,x\to y)$ in $\treespr$, so $h\in\cH(\trees)$, and in fact by the consistency conditions $h\in\cH_\star(\trees)$. Thus there is a bijection between BP fixed points $h\in\cH_\star(\trees)$ and translation-invariant Markovian $\nu^h\in\gibbs_\trees$. In particular, for two-spin models, the $\nu^\pm\in\gibbs_\trees$ of \S\ref{ss:intro.loc} and the $h^\pm\in\cH_\star$ of \S\ref{ss:bethe} are related by this correspondence and so may be regarded as essentially equivalent.

The main implication of Lem.~\ref{l:int} is the following (which may also be obtained as a special case of \cite[Thm.~1.13]{\DMS}): if for $B\in[B_0,B_1]$ we have $h\equiv h(B)\in\cH^B_\star$ which is continuous and of bounded total variation in $B$, then
\beq\label{e:a.ubd}
\limsup_{n\to\infty} a^\vertex_n(B) \le a^\vertex_\trees(B,h),
\eeq
implies $\limsup_{n\to\infty}[\phi_n(B_1)-\phi_n(B_0)]\le\Phi_\trees(B_1)-\Phi_\trees(B_0)$. (The statement also holds with $B,a^\vertex$ replaced by $\be,a^\edge$.) From the above discussion we can re-express
$$a^\vertex_\trees(B,h)\equiv
	\E_\trees[\angl{ \pd_B\vxi(\si_\rt)}^B_{\nu^h}],\quad
a^\edge_\trees(\be,h)\equiv
	\f12\E_\trees\Big[ \sum_{j\in\pd\rt}
		\angl{\pd_\be\xi(\si_\rt,\si_j)}^\be_{\nu^h} \Big],
$$
so showing \eqref{e:a.ubd} amounts to proving a relation between the expectation of a local observable, in this case $\pd_B\vxi(\si_i)$, in the finite graph to the expectation of the analogous observable $\pd_B\vxi(\si_\rt)$ under translation-invariant Markov chains on the limiting tree. In the remainder of the section we carry out this scheme. Note that for the Ising model
$$a^\vertex_n(B)=\E_n[\angl{\si_{I_n}}^B_n],\quad
a^\edge_n(\be)=\f12\E_n\Big[\sum_{j\in\pd I_n}\angl{\si_{I_n}\si_j}^\be_n\Big].$$
For the hard-core model $a^\vertex_n(B)=\E_n[\angl{\bar\si_{I_n}}^B_n]$.

\subsubsection{Interpolation for hard-core}
\label{sss:is}

\blem\label{l:hc.opt}
For the hard-core model at fugacity $\lm$, the supremum of $\angl{\bar\si_\rt + d^{-1}\sum_{j\in\pd\rt}\bar\si_j}_\mu$ over $\mu\in\gibbs_\treereg$ is achieved precisely by the measures $\mu^\pm$. \\
Consequently, the supremum of $\bipE[\angl{\bar\si_\rt}_\nu]$ over translation-invariant $\nu\in\gibbs_\trees$ is achieved precisely by the $\nu^\pm$.
\bpf
By extremal decomposition, assume without loss that $\mu$ is itself extremal, with $h^\mu\equiv(h^\mu_{x\to y})$ as defined above. For $j\in\pd\rt$ write $q_j\equiv h^\mu_{j\to\rt}(-)$: then
\begin{align*}
\Big\langle
	\bar\si_\rt + \f{1}{d}
	\sum_{j\in\pd\rt} \bar\si_j\Big\rangle_{\mu}
&= \f{\lm \prod_{j\in\pd\rt} q_j
	+ d^{-1} \sum_{j\in\pd\rt} (1-q_j)}
	{1 + \lm \prod_{j\in\pd\rt} q_j}\\
&= 1 - \f{1}{d} \f{\sum_{j\in\pd\rt} q_j}{1 + \lm \prod_{j\in\pd\rt} q_j}.
\end{align*}
For fixed $\sum_j q_j$ this is (strictly) maximized by taking all $q_j\equiv q$, so the above is
$$\le 1-\f{1}{\max_{q^-\le q\le q^+}[q^{-1}+\lm q^{d-1}]},$$
where $q^-$ and $q^+$ are the minimal and maximal values for $q$, corresponding to $\mu^-$ and $\mu^+$ respectively. Since $q^{-1}+\lm q^{d-1}$ is convex, the maximum can only be attained at the endpoints, and in fact it is attained at both endpoints with value $(1/q^-+1/q^+-1)^{-1}$.\epf
\elem

\blem\label{l:hc.cty}
For the hard-core model, $\cH^\lm_\star(\set{\treereg})$ consists of a single message $h^\star\equiv h^\star(\lm)$. $\cH^\lm_\star(\set{\treereg_\pm})$ consists of the messages $h^\star,h^\pm$ which coincide for $\lm\le\lm_c$ and are distinct for $\lm>\lm_c$. The messages are continuous in $\lm$, smooth except possibly at $\lm=\lm_c$.

\bpf
For the hard-core model, the function $\BPF(h)$ ($h\in\simplex$) of \eqref{e:bpf} is expressed in terms of $q\equiv h(-)$ as $\BPF(q) = (1+\lm q^{d-1})^{-1}$. As $q$ increases from $0$ to $1$, $\BPF$ decreases from $1$ to $(1+\lm)^{-1}$, so $\BPF$ has a unique fixed point $q_\star$ which is smoothly decreasing in $\lm$. We compute
$$\BPF' = -\f{(d-1)}{q} \BPF [1-\BPF],\quad
\BPF''
=- \f{d-1}{q^2} \BPF (1-\BPF)
	[(d-1)2\BPF-d],$$
so $\BPG\equiv\BPF^{(2)}$ has second derivative
\begin{align*}
\BPG''
&= (\BPF'\circ\BPF)  \BPF'' + (\BPF''\circ\BPF)  (\BPF')^2
= \frac{(d-1)^2}{q^2} (1-\BPF)(1-\BPG) \BPG \cdot Q,\\
Q&\equiv -2(d-1)^2(1-\BPF)\BPG
	+(d-2) [d(1-\BPF)+\BPF].
\end{align*}
Setting this to zero gives
$$\BPG = \f{(d-2)}{2(d-1)^2}
	\lp d + \f{\BPF}{1-\BPF} \rp.$$
The left-hand side is decreasing in $\BPF$ while the right-hand side is increasing, so $\BPG$ has at most one inflection point, hence at most three fixed points. If $\BPG$ has a fixed point which is not equal to $q_\star$, then necessarily it has exactly three fixed points $\BPF(q_\circ)<q_\star<q_\circ$, so $\BPG'(q_\star)>1$. But
$$\BPG'(q_\star)=\BPF'(q_\star)^2 = (d-1)^2(1-q_\star)^2$$
is smoothly increasing in $\lm$ with $\BPG'(q_\star)=1$ precisely at $\lm=\lm_c(d)$, so we see $\BPG$ has a unique fixed point $h_\star$ when $\lm\le\lm_c$, and when $\lm>\lm_c$ it has three fixed points $\BPF(q_\circ)<q_\star<q_\circ$ which are smooth on the open interval $(\lm_c,\infty)$.

It remains to verify that $q_\circ\to q_\star$ as $\lm\decto\lm_c$. Suppose otherwise, so that
$$\limsup_{\lm\decto\lm_c} \BPF(q_\circ)+2\ep
	<q_\star
	<\liminf_{\lm\decto\lm_c} q_\circ-2\ep$$
for some $\ep>0$. It is possible to take a sequence $\lm\decto\lm_c$ along which the inflection point of $\BPG$ always lies on the same side of $q_\star$: assume it is $\le q_\star$ (the argument for the $\ge q_\star$ case is symmetric), so that $\BPG'$ is decreasing on $q\ge q_\star$. By the mean value theorem applied to the interval $[q_\star+\ep,q_\star+2\ep]$,
$$\BPG'(q_\star+\ep)
\ge \f{\BPG(q_\star+2\ep)-\BPG(q_\star+\ep)}{\ep}
\ge \f{q_\star+2\ep-[q_\star+\BPG'(q_\star)\ep]}{\ep}
= 2-\BPG'(q_\star).$$
Therefore $2-\BPG'(q_\star)\le\BPG'(q)\le\BPG'(q_\star)$ for all $q\in[q_\star,q_\star+\ep]$, and consequently $\BPF'(q)=1$ for all $q\in[q_\star,q_\star+\ep]$ at $\lm_c$. This gives the desired contradiction and the lemma follows.
\epf
\elem

\bppn\label{p:hc.z}
For the hard-core model,
\bnm[(a)]
\item The definitions \eqref{e:Bethe} and \eqref{e:bethe.af} of $\Phi\equiv\Phi_{\set{\treereg_\pm}}$ coincide. \\ If $G_n\lpc\treereg$ then $\phi=\Phi$ for $\lm\le\lm_c$ and $\limsup_n\phi_n\le\Phi$ for $\lm>\lm_c$.
\item If $G_n\lwc\biptreereg$ then $\phi=\Phi$ for all $\lm>0$.
\enm

\bpf
(a) Since any subsequential local weak limit of the measures $\nu_n$ must be translation-invariant, the second part of Lem.~\ref{l:hc.opt} implies
$$\limsup_{n\to\infty}
	a^\vertex_n(B)
=\limsup_{n\to\infty} \E_n[\angl{ \bar\si_{I_n}}^B_n]\\
\le \bipE[\angl{\bar\si_\rt}_{\nu^\pm}]
	\equiv a^\vertex(B,h^\pm)$$
(the inequality can alternatively be obtained by expressing $a^\vertex_n(B)$ as $o(1)+\E_n[ \angl{
	\bar\si_{I_n}
	+ d^{-1}\sum_{j\in\pd I_n}\bar\si_j}^B_n ]/2$ and directly applying the first part of Lem.~\ref{l:hc.opt}). Also, $a^\vertex(B,h^\star)\le a^\vertex(B,h^+)=a^\vertex(B,h^-)$, with equality for $\lm\le\lm_c$ and with strict inequality for $\lm>\lm_c$. It then follows from Lem.~\ref{l:int} (using Lem.~\ref{l:hc.cty}) that for $\lm_c<\lm$,
$$\limsup_n[\phi_n(\lm)-\phi_n(\lm_c)]
\le\Phi(\lm,h^\pm)-\Phi(\lm_c,h^\pm)
> \Phi(\lm,h^\star)-\Phi(\lm_c,h^\star).$$
It was shown in \cite[Thm.~1.11]{\DMS} that $\phi=\Phi(h^\star)=\Phi(h^\pm)$ for $\lm\le\lm_c$ so the claim follows (again making use of Lem.~\ref{l:hc.cty}).

\medskip\noindent
(b) It suffices to show that for $G_n\lwc\biptreereg$,
\beq\label{e:hc.infty}
\lim_{\lm\to\infty}\liminf_{n\to\infty}(\phi_n-\Phi)\ge0.
\eeq
Indeed,
\begin{align*}
\Phi^\vertex(\lm) &= \f{1}{2}\log[\lm (q^+)^d+1] + \f{1}{2}\log[\lm (q^-)^d+1],\\
\Phi^\edge(\lm) &= \f{d}{2}\log(1-(1-q^-)(1-q^+)),
\end{align*}
and $\lim_{\lm\to\infty}q^+=\lim_{\lm\to\infty}(1-q^-)=1$, so $\lim_{\lm\to\infty}\Phi(\lm)-[\log(\lm+1)]/2=0$. But
$$\liminf_{n\to\infty}\phi_n
\ge \liminf_{n\to\infty}
\f{1}{n} \E_n \bigg[
	\log
	\sum_{j=0}^{\al_n} {n\choose j} \lm^j \bigg]
= \log(1+\lm) \f{\E_n[\al_n]}{n}$$
for $\al_n$ the independence number of $G_n$. But $\al_n$ is at least the number of black vertices with no black neighbors, so $G_n\lwc\biptreereg$ implies $\limsup_n\E_n[\al_n]/n\le1/2$. This proves \eqref{e:hc.infty} from which the result follows.
\epf
\eppn

\subsubsection{Interpolation for Ising}

\blem\label{l:ising.v.opt}
For the Ising model with parameters $\be<0$ and $B\in\R$, the supremum of $\angl{\si_\rt + d^{-1} \sum_{j\in\pd\rt}\si_j}_\mu$ over $\mu\in\gibbs_\treereg$ is achieved precisely by the $\mu^\pm$. \\
Consequently, the supremum of $\bipE[\angl{\si_\rt}_\nu]$ over translation-invariant $\nu\in\gibbs_\trees$ is achieved precisely by the $\nu^\pm$.

\bpf
We argue as in the proof of Lem.~\ref{l:hc.opt}: assume $\mu$ is extremal, and write $(h_j,1-h_j)\equiv(h^\mu_{j\to\rt}(+),h^\mu_{j\to\rt}(-))$: then
$$\Big\langle
\si_\rt + \f{1}{d} \sum_{j\in\pd\rt}\si_j
\Big\rangle_\mu
= \f{[e^B - e^{-B} \prod_{k\in\pd\rt}R_k]
	+ d^{-1}\sum_{j\in\pd\rt}
		[e^B A_j + (e^{-B} \prod_{k\in\pd\rt} R_k) B_j]}
	{ e^B + e^{-B} \prod_{k\in\pd\rt} R_k}$$
where $R_j\equiv [e^{-\be}h_j+e^\be(1-h_j)]/[e^\be h_j + e^{-\be}(1-h_j)]$, and
\begin{align*}
A_j \equiv \f{e^\be h_j-e^{-\be}(1-h_j)}{e^\be h_j+e^{-\be}(1-h_j)}
&= \f{2 R_j - (e^{2\be}+e^{-2\be})}{e^{-2\be}-e^{2\be}},\\
B_j\equiv \f{e^{-\be}h_j-e^\be(1-h_j)}{e^{-\be}h_j+e^\be(1-h_j)}
&= \f{(e^{2\be}+e^{-2\be})-2R_j^{-1}}{e^{-2\be}-e^{2\be}}.
\end{align*}
If $(h_k)_{k\ne j}$ are fixed, the expression is maximized by taking $h_j$ as large or as small as possible. Since $\be<0$, for any fixed $\prod_{k\in\pd\rt} R_k$, both $\sum_{j\in\pd\rt} A_j$ and $\sum_{j\in\pd\rt} B_j$ are maximized by taking all the $R_j$ equal, i.e.\ with $h_j\equiv h$. The overall maximum is then attained for $h$ equal to $h^+$ or $h^-$, the minimal and maximal values for $q$ corresponding to the $\mu^\pm$. In fact it is attained at both endpoints with value
$$\f{2(h^-+h^+-1)}
	{1 + (e^{-2\be}-1)(h^-+h^+)
		-2(e^{-2\be}-1)h^- h^+},$$
which concludes the proof.
\epf
\elem

\blem\label{l:ising.e.opt}
For the Ising model with parameters $\be<0$ and $B=0$, the supremum of $-\angl{\sum_{j\in\pd\rt}\si_\rt\si_j}_\mu$ over $\mu\in\gibbs_\treereg$ is achieved precisely by the $\mu^\pm$.

\bpf
As in the proof of Lem.~\ref{l:ising.v.opt}, assume $\mu$ is extremal and write $(h_{x\to y},1-h_{x\to y})\equiv( h^\mu_{x\to y}(+),h^\mu_{x\to y}(-))$. Then
\beq\label{e:ising.e}
-\angl{\si_\rt\si_j}_\mu
= \f{e^{-\be}
	-(e^\be+e^{-\be})[h_{\rt\to j} h_{j\to\rt}+(1-h_{\rt\to j})(1-h_{j\to\rt})]}
	{e^{-\be}-(e^{-\be}-e^\be)
		[h_{\rt\to j} h_{j\to\rt}+(1-h_{\rt\to j})(1-h_{j\to\rt})]}.
\eeq
The partial derivative with respect to $h_{j\to\rt}$ has the same sign as $1/2-h_{\rt\to j}$, so $-\angl{\si_\rt\si_j}_\mu$ is maximized with $(h_{\rt\to j},h_{\rt\to j})$ equal to $(h^-,h^+)$ and $(h^+,h^-)$, corresponding to the measures $\mu^\pm$.
\epf
\elem

\blem\label{l:ising.cty}
For the anti-ferromagnetic Ising model, $\cH^{\be,B}_\star(\set{\treereg})$ consists of a single message $h^\star\equiv h^\star(\be,B)$. $\cH^\lm_\star(\set{\treereg_\pm})$ consists of the messages $h^\star,h^\pm$ which coincide for $\be_{c,\af}(B,d)\le\be\le0$ and are distinct for $\be<\be_{c,\af}(B,d)$. The messages are continuous in $\be,B$, smooth except possibly where $\be=\be_{c,\af}(B,d)$.

\bpf
For the Ising model, the function $\BPF(h)$ ($h\in\simplex$) of \eqref{e:bpf} is expressed in terms of $t\equiv \log[h(+)/h(-)]$ as
\beq\label{e:ising.bp}
\BPF(t) = 2B+(d-1)\log\ls\f{e^t + \thet}{ \thet e^t + 1}\rs,\quad \thet\equiv e^{-2\be}.
\eeq
$\BPF-2B$ is an odd function of $t\in\R$, identically zero when $\be=0$ and strictly monotone otherwise, going from $(d-1)\log\thet$ to $-(d-1)\log\thet$ as $t$ increases from $-\infty$ to $\infty$. Suppose from now on that $\be<0$. Then $\BPF$ has a unique fixed point $t_\star$ of the same sign as $B$, smoothly increasing in $B$, and smooth in $\be$ with absolute value increasing as $\be$ becomes more negative.

For $\BPG\equiv\BPF^{(2)}$ we compute
\begin{align*}
\BPG''(t)&=
\frac{(d-1)^2 \thet (\thet^2-1)^2 e^{\BPF+t}
	(A_+ e^{2\BPF} + A e^\BPF + A_-)}
	{(e^\BPF+\thet)^2
	(\thet e^\BPF +1)^2
	(e^t + \thet)^2
	(\thet e^t+1)^2},\\
A_\pm &\equiv
	-\thet(e^{2t}-1) \pm (d-1)(\thet^2-1) e^t,\\
A &\equiv -(\thet^2+1)(e^{2t}-1).
\end{align*}
If $A_+=0$ then clearly $\BPG$ can have at most one inflection point, so suppose $A_+>0$: then setting $\BPG''$ to zero results in $e^\BPF = r_\pm$ where
\begin{align*}
r_\pm &\equiv (-A \pm \sqrt{D})/(2A_+),\\
D &\equiv A^2-4 A_+ A_- 
= (\thet^2-1)^2
	[1 +4 ((d - 1)^2 - 1/2) e^{2t} + e^{4t}]>0.
\end{align*}
If $A/A_+\ge0$ then at most one of the $r_\pm$ can be positive. If $A\ge0$ then $t\le0$ which implies $A_+>0$, so it remains only to consider the case $A<0<A_+$: in this case, $A<0$ implies $t>0$ and so $A_-<0$, therefore
$$\abs{D}-A^2=-4 A_+ A_- > 0$$
which implies $r_-<0<r_+$. Thus $\BPG$ has at most one inflection point for any $\be<0$.

By implicit differentiation we find
\begin{align*}
\pd_\thet [\BPF'(t_\star)]
&= \left.[ \pd_\thet \BPF'(t)]\right|_{t=t_\star}
	+ \BPF''(t_\star) \pd_\thet t_\star
= \left.[ \pd_\thet \BPF'(t)]\right|_{t=t_\star}
	+ \BPF''(t_\star)
	\f{\left.[ \pd_\thet \BPF(t)]\right|_{t=t_\star}}
	{1-\BPF'(t_\star)}\\
&=- \left.\f{(d-1) e^t}{(e^t+\thet)(\thet e^t+1)}
	\f{[d(\thet^2-1)+2](e^{2t}+1)+4\thet e^t }
	{[d(\thet^2-1)+2]e^t+\thet(e^{2t}+1)}\right|_{t=t_\star}
<0,
\end{align*}
so $\BPG'(t_\star)=\BPF'(t_\star)^2$ increases smoothly as $\be$ decreases. The result then follows by repeating the argument of Lem.~\ref{l:hc.cty}.
\epf
\elem

\bppn\label{p:ising.z.af}
For the anti-ferromagnetic Ising model,
\bnm[(a)]
\item The definitions \eqref{e:Bethe} and \eqref{e:bethe.af} of $\Phi\equiv\Phi_{\set{\treereg_\pm}}$ coincide. \\ If $G_n\lpc\treereg$ then $\phi=\Phi$ for $\be\ge\be_{c,\af}(B)$ and $\limsup_n\phi_n \le\Phi$ for $\be<\be_{c,\af}(B)$.
\item If $G_n\lwc\biptreereg$ then $\phi=\Phi$ for all $\be,B$.
\enm

\bpf
(a) First fix $B=0$: since any subsequential local weak limit of the measures $\nu_n$ must be translation-invariant, Lem.~\ref{l:ising.e.opt} gives
\begin{align*}
&\liminf_{n\to\infty}
a^\edge_n(\be)
=\f12 \liminf_{n\to\infty}
	\E_n\Big[
	\sum_{j\in\pd I_n} \angl{\si_{I_n}\si_j}^\be_n
	\Big]\\
&\ge \f{1}{2}\bipE\Big[ \sum_{j\in\pd\rt} \angl{\si_\rt\si_j}_{\nu^\pm} \Big]
\ge a^\edge(\be,h^\pm).
\end{align*}
Also, $a^\edge(\be,h^\star)\ge a^\edge(\be,h^\pm)$ with equality for $\be\ge\be_{c,\af}$ and with strict inequality for $\be<\be_{c,\af}$. It then follows from Lem.~\ref{l:int} (together with Lem.~\ref{l:ising.cty} and the previous result for $\be\ge0$) that for $B=0$ and $\be\le0$ we have
\beq\label{e:ising.ineq}
\limsup_{n\to\infty}\phi_n \le \Phi(h^\pm)>\Phi(h^\star).
\eeq
Using Lem.~\ref{l:ising.v.opt} to interpolate in $B$ (as in the proof of Propn.~\ref{p:hc.z}) then gives \eqref{e:ising.ineq} for all $\be\le0, B\in\R$.

\medskip\noindent
(b) Consider the limits $\be\to-\infty$ and $B\to\infty$:
$$\lim_{\be\to-\infty}[\Phi(\be,0)+\be d/2]=0,\quad
\lim_{B\to\infty}[\Phi(\be,B)-B-\be d/2]=0.$$
If $G_n\lwc\biptreereg$, then
\begin{align*}
\liminf_{n\to\infty} \phi_n
&\ge \liminf_{n\to\infty} -\be\f{\E_n[\abs{E_n}]}{n} = -\be d/2,\\
\liminf_{n\to\infty} \phi_n
&\ge B+ \liminf_{n\to\infty} \be\f{\E_n[\abs{E_n}]}{n} = B + \be d/2,
\end{align*}
so in fact $\phi=\Phi$ for all $\be,B\in\R$.
\epf
\eppn

For completeness we review what is known for the ferromagnetic Ising model:

\bppn\label{p:ising.z.fer}
For the ferromagnetic Ising model on $G_n\lpc\treereg$, $\phi$ exists and equals $\Phi_{\set{\treereg}}$ as defined by \eqref{e:Bethe} (and given more explicitly by \eqref{e:bethe.f}).

\bpf
In this setting $\cH_\star$ corresponds simply to the fixed points of a single iteration of the map $\BPF$ of \eqref{e:ising.bp}, which is analyzed for example in \cite[Lem.~4.6]{\DMS}. $\BPF$ always has between one and three fixed points, and we write $t^+$ and $t^-$ for the maximal and minimal fixed points respectively, corresponding to $h^+$ and $h^-$ as in \eqref{e:bethe.f}. By symmetry we may always suppose $B\ge0$.

At $B=0$, a fixed point is always given by $t^\free=0$, unique provided $\BPF'(0)\le1$. $\BPF'(0)$ increases monotonically in $\be$ and reaches $1$ at $\be_{c,\fer}(B=0,d)$, and for $\be>\be_{c,\fer}(0,d)$ there are three distinct fixed points $t^-<0=t^\free<t^+$, with $t^+ = -t^- \decto 0=t^\free$ as $\be\decto\be_{c,\fer}(0,d)$. Since adding $B$ simply shifts the map $\BPF$ of \eqref{e:ising.bp} by the constant $2B$, it is easy to deduce the behavior for general $B\ge0$: if $\BPF'(t)|_{t=0}\le1$ then $\BPF$ has a unique fixed point $t^+=t^-$ which is zero when $B=0$ and increases smoothly in $B$. If $\BPF'(t)|_{t=0}>1$, then at $B=0$ the map $\BPF$ has three fixed points $t^-<t^\free=0<t^+$. As $B$ increases, $t^\pm$ increase smoothly while $t^\free$ decreases smoothly. The fixed points $t^-$ and $t^\free$ merge at the threshold $B=B_{c,\fer}(\be,d)$, and for $B$ above this threshold we again have $t^-=t^+$.

It follows from \cite[Thm.~2.4]{\DemboMontanariIsing} (see also \cite[Thm.~1.8]{\DMS}) that for $G_n\lpc\treereg$, $\phi$ exists and equals $\Phi$ as defined by \eqref{e:bethe.f}. Therefore it remains to verify that
\beq\label{e:ising.fer}
\Phi(h^+)\ge\Phi(h^-)\vee\Phi(h^\free)
\quad\text{for all }B\ge0.
\eeq
At $B=0$, it follows from the above that $\Phi(h^\free)$ and $\Phi(h^+)=\Phi(h^-)$ are continuous in $\be$. Writing $h\equiv h(+)$ and recalling \eqref{e:ising.e}, we compute
$$\pd_h \angl{\si_\rt \si_j}
= \f{4 e^{2\be} (2 h-1)}
	{[e^{2\be}- 2(e^{2\be}-1) h(1-h) ]^2},$$
so $\angl{\si_\rt \si_j}_{h^+}\ge\angl{\si_\rt \si_j}_{h^\free}$, hence $a^\edge(\be,h^\pm)\ge a^\edge(\be,h^\free)$, and then Lem.~\ref{l:int} gives $\Phi(h^\pm)\ge\Phi(h^\free)$ at $B=0$. Next, for all $B\ge0$, clearly $a^\edge(\be,h^+) \ge a^\edge(\be,h^-)\vee a^\edge(\be,h^\free)$, so another application of Lem.~\ref{l:int}  gives \eqref{e:ising.fer} from which the proposition follows.
\epf
\eppn

\bpf[Proof of Thms.~\ref{t:intro.z} and~\ref{t:z}]
Follows by combining the reduction of \S\ref{ss:reduc.two} with the results of Propns.~\ref{p:hc.z},~\ref{p:ising.z.af}, and~\ref{p:ising.z.fer}.
\epf

\section{Local structure of measures}\label{s:loc}

In this section we show how Thm.~\ref{t:z} can be used to deduce Thm.~\ref{t:loc} by straightforward modifications of the arguments of~\cite{\MMS}.

\bpf[Proof of Thm.~\ref{t:loc}~(a)]
Observe that $\pd_B^2\phi_n = n^{-1}\E_n[ \angl{S^2}^B_n-(\angl{S}^B_n)^2 ]$ where $S$ is $\sum_{i\in V_n}\si_i$ for Ising and $\sum_{i\in V_n}\bar\si_i$ for hard-core, so the $\phi_n$ are convex and hence so is the limit $\phi$. Convex functions are absolutely continuous, so it holds for a.e.\ $B$ that $\phi_n,\phi$ are differentiable in $B$ with $\pd_B\phi_n\to\pd_B\phi=\pd_B\Phi$ (by Thm.~\ref{t:z}). It follows from \cite[Propn.~2.4]{\DMS} that $\pd_B\Phi \equiv a^\vertex(B,h^+) = \bipE[\angl{\si_\rt}_{\nu^+}]$. But for any subsequential local weak limit $\nu$ of the $\nu_n$, we also have
$$\pd_B\phi_n \equiv a^\vertex_n(B)
	=\E_n[\angl{\si_{I_n}}^B_n]
\to \bipE[ \angl{\si_\rt}_\nu].$$
Therefore $\bipE[ \angl{\si_\rt}_\nu]=\bipE[ \angl{\si_\rt}_{\nu^+}]$, and it follows from Lem.~\ref{l:hc.opt} and Lem.~\ref{l:ising.v.opt} that $\nu$ is a convex combination of the $\nu^\pm$. Since the $G_n$ are symmetric, we must have $\nu=(\nu^++\nu^-)/2$.
\epf

We now analyze the conditional measures $\nu^\pm_n$, beginning with an easy observation:

\blem\label{l:zero}
For anti-ferromagnetic two-spin models on $G_n\lwc\biptreereg$,
$$\lim_{n\to\infty}
\E_n\Big[\nu_n\Big(\sum_{i\in V_n}\col_i\si_i=0\Big)\Big]=0.$$

\bpf
For the Ising model see \cite[Lem.~4.1]{\MMS}. For the hard-core model, let $A_n$ denote the set of vertices $i\in V_n$ with $\ball{2}{i}$ isomorphic to $\treereg_+^2$, the depth-two subtree of $\treereg_+$; then $A_n$ is necessarily an independent set of black vertices. The probability that $\sum_{i\in A_n}\col_i\bar\si_i=\sum_{i\in A_n}\bar\si_i$ takes value $j$, conditioned on all the spins $(\bar\si_i)_{i\notin A_n}$, is $\P(X=j)$ where $X$ is a binomial random variable on $N=\abs{\set{i\in A_n : \bar\si_{\pd i}\equiv0}}$ number of trials with success probability $\lm/(1+\lm)$. If $N\ge\ep n$ then $\P(X=j)=O(1/\sqrt{\ep n})$ uniformly in $j$ (e.g.\ by the Berry-Ess\'een theorem). If $N<\ep n$ then $\sum_{i\in\pd A_n} \bar\si_i \ge (\abs{A_n}-\ep n)/d$, so
$$\f1n\sum_{i\in V_n}\col_i\si_i
=\f2n\sum_{i\in V_n}\col_i\bar\si_i
	-\f1n\sum_{i\in V_n}\col_i
<\ep - \f{\abs{A_n}/n-\ep}{d}
	+\f{\abs{V_n \setminus (A_n\cup\pd A_n)}}{n}
	-\f1n\sum_{i\in V_n}\col_i.$$
As $n\to\infty$ the right-hand side tends in probability to $[-1/2+\ep(d+1)]/d$, which is negative for small $\ep$. Combining the above observations concludes the proof for the hard-core model.
\epf
\elem

In view of Lem.~\ref{l:zero} we may without loss restrict attention to the measures $\nu^+_n$. Define the local functions (cf.~\cite[eq.~(3.9)]{\MMS})
$$F^t_i \equiv F^t_i(\de,\vec\si)
\equiv \I \Big\{
	\sum_{j\in\ball{t}{i}} \col_j\si_j
	\le-\de\abs{\ball{t}{i}}
	\Big\};$$
$F^t_i$ indicates the vertices of $G_n$ which are locally not in the $+$ phase.

\bpf[Proof of Thm.~\ref{t:loc}~(b)]
We outline the steps of the proof of \eqref{e:locConvergence} following~\cite{\MMS}, describing modifications where needed.
\bitm
\item Let $\nu^*$ denote any subsequential local weak limit of the $\nu^+_n$. Then $\nu^*\in\gibbs_\trees$ (see \cite[Lem.~3.4]{\MMS}). By Lem.~\ref{l:zero}, $\nu^+_n$ has free energy density converging to $\phi$, so the proof of Thm.~\ref{t:loc}~(a) implies that $\nu^*=(1-q)\nu^++q \nu^-$ for some $q\in[0,1]$.

\item By local weak convergence, $\lim_{n\to\infty} \E_n[ \angl{F^t_{I_n}}_n ] = \bipE[\angl{F^t_\rt}_{\nu^*}]$; further, if $J_n$ denotes a uniformly random neighbor of $I_n$, then
$$\lim_{n\to\infty} \E_n[ \angl{\Ind{F^t_{I_n}\ne F^t_{J_n}}}_n ]
= \bipE[\angl{\Ind{F^t_\rt\ne F^t_j}}_{\nu^*}],\quad j\in\pd\rt$$
(cf.~\cite[Lem.~3.7]{\MMS}).

\item For the hard-core or anti-ferromagnetic Ising model in non-uniqueness regimes, there exists $\de>0$ such that
\begin{align*}
&\lim_{t\to\infty} \angl{F^t_\rt}_{\nu^+} = 0
= 1-\lim_{t\to\infty} \angl{F^t_\rt}_{\nu^-},\\
&\lim_{t\to\infty} \angl{F^t_\rt\ne F^t_j}_{\nu^+}
= 0 = \lim_{t\to\infty} \angl{F^t_\rt\ne F^t_j}_{\nu^-}
\end{align*}
(cf.~\cite[Lem.~3.8]{\MMS}).
It follows that for sufficiently large $t$
$$\lim_{n\to\infty} \E_n[ \angl{F^t_{I_n}}_n ]
	\ge q-\ep,\quad
\lim_{n\to\infty} \E_n[ \angl{\Ind{F^t_\rt\ne F^t_j}}_n ]\le\ep.$$
\eitm
The argument of \cite[Propn.~3.9]{\MMS} (using the edge-expansion hypothesis) now gives a contradiction unless $q=0$ establishing~\eqref{e:locConvergence}.  The proof of \eqref{e:magnetization} then follows from applying the proof of \cite[Thm.~2.5]{\MMS} to the bipartite case.
\epf

\section{Computational hardness}\label{s:comp}

In this section we construct the bipartite expander gadgets to be used in the reduction to $\maxcut$ (Lem.~\ref{l:exp}) and refine Thm.~\ref{t:loc} to an approximate conditional independence statement for the gadgets (Propn.~\ref{p:product}). We conclude with the proof of our main results Thms.~\ref{t:hc.comp} and~\ref{t:ising.comp}.

For any fixed positive integer $k$, $G^k_{2n}$ will be a bipartite graph on $2n$ vertices with $n$ even, defined as follows:
\bitm
\item Let $H_n$ be a graph on $n$ vertices of maximum degree $d$, generated by the configuration model as follows: take a uniformly random matching $\match$ of $[dn]$, and put an edge $(i,j)$ in $H_n$ for every edge $(i',j')\in\match$ with $i'\in i+n\Z$, $j'\in j+n\Z$ (self-loops and multi-edges allowed).
\item Take $G_{2n}$ to be the bipartite double cover of $H_n$: the two parts of $G_{2n}$ are $(i_+)_{i=1}^n$ and $(i_-)_{i=1}^n$, and we put two edges $(i_+,j_-)$ and $(j_+,i_-)$ in $G_{2n}$ for every edge $(i,j)\in H_n$ (multi-edges allowed).
\item Choose $k$ vertices $(i^\ell)_{\ell=1}^k$ uniformly at random from $H_n$, and for each $\ell$ choose $j^\ell\in\pd i^\ell$ uniformly at random. $G^k_{2n}$ is the simple bipartite graph formed by deleting the edges $(i^\ell_\pm,j^\ell_\mp)$ from $G_{2n}$ and merging any remaining multi-edges in the graph into single edges. Write $W^\pm\equiv\set{i^\ell_\pm,j^\ell_\pm}_{\ell=1}^k$ and $W\equiv W^+\cup W^-$.
\eitm
The graphs $G_{2n}$ are $d$-regular with probability bounded away from zero as $n\to\infty$ (see e.g.\ \cite[Ch.~9]{\JansonEtAl}).

\blem\label{l:exp}
Let $k$ be fixed. For all $\de>0$ there exists $\lm_\de>0$ such that the $G^k_{2n}$ are $(\de,1/2,\lm_\de)$-edge expanders with high probability as $n\to\infty$.

\bpf
By stochastic domination we may assume $d=3$. For $S\subset H_n$ with $\abs{S}=m$, the probability that there are exactly $j$ edges in $H_n$ between $S$ and its complement is
$$P_{j,m}
= I_{j,m}
\f{{3m\choose j}{3(n-m)\choose j} j!
	M_{3(m-j)} M_{3(n-m-j)}}{M_{3n}},$$
where $I_{j,m}$ is the indicator that $m-j$ is even, and $M_\ell=(\ell-1)!!
=\pi^{-1/2}\Gam[(\ell+1)/2] 2^{\ell/2}$ is the number of matchings on $[\ell]$ for $\ell$ even. By Stirling's approximation, if $\de\le m/n\le 1-\de$ and $j=\gam n$, then
$$P_{j,m} = I_{j,m}\exp\Big\{ -n \Big[
	\f32 H(m/n) -\gam\log\gam +
		O_\de(\gam)
	\Big]+o_\de(n) \Big\}$$
(where $H(p)$ denotes the binary entropy function $-p\log p-(1-p)\log(1-p)$). There are $\le e^{n H(m/n)}$ subsets of $H_n$ of size $m$ so there exists $\gam_\de>0$ such that with probability at least $n e^{-n H(\de)/4}$, all subsets of $H_n$ of size between $\de n$ and $(1-\de)n$ have expansion at least $\gam_\de$.

We now show expansion for $G^k_{2n}$: since $k$ does not change with $n$ and the number of edges leaving any set of vertices decreases by at most a factor of $3$ when multi-edges are merged into single edges, it suffices to show expansion for $G_{2n}$. Let $S_\pm$ be subsets of the $\pm$ sides of $G_{2n}$ such that $S\equiv S_+\cup S_-$ has size $\le n$. If the projection $\pi S$ of $S$ in $H_n$ has size $\le(1-\de)n$, then $S$ has expansion at least $\gam_\de/2$. Suppose $\abs{\pi S}\ge(1-\de)n$: without loss $\abs{S_+}\ge\abs{S_-}$, so $\abs{\pi S_+\setminus \pi S_-}\ge(1/2-\de)n$. If there are fewer than $\gam\abs{S}$ edges leaving $S$, then there must be at least $3(1/2-\de)n-\gam n$ edges between $\pi S_+\setminus \pi S_-$ and its complement in $H_n$. A similar analysis as above shows that for sufficiently small $\de$ there exists $\gam_\de>0$ such that the probability $G_{2n}$ has such a set $S$ is $\le e^{-n(\log 2)/4}$, and this concludes the proof.
\epf
\elem

Recall that we use $W^\pm$ to denote the endpoints on the $\pm$ sides of the $2k$ edges deleted from $G_{2n}$ in the formation of $G^k_{2n}$. Recall also the definitions of $\mu^\pm\in\gibbs_\treereg$, and write $h^\pm\equiv h^{\mu^\pm}_{\rt\to j}\in\simplex$. For $h,h'\in\simplex$ define $h\otimes_\psi h'\in\simplex_{\spins^2}$ by
\beq\label{e:z.psi.h}
(h\otimes_\psi h')(\si,\si')
= \f{h(\si)\psi(\si,\si')h(\si')}{z(h\otimes_\psi h')},
\eeq
for $z(h\otimes_\psi h')$ the normalizing constant.

\bppn\label{p:product}
The conditional measure $\nu_{G^k_{2n}}^\pm(\vec\si_W=\cdot)$ converges to the product measure
$$Q^\pm_W(\vec\si)
\equiv
\prod_{w\in W^+} h^\pm(\si_w)
\prod_{w\in W^-} h^\mp(\si_w).$$

\bpf
Let $B_t$ denote the union of the balls $\ball{t}{w}\subseteq G_{2n}$ over $w\in\set{i^\ell_\pm}_{\ell=1}^k$; assume that $B_t$ is a disjoint union of graphs isomorphic to $\treereg^t$ with internal boundary $S_t\equiv B_t\setminus B_{t-1}$, which is the case with high probability. For $\vec\eta\in\spins^{S_t}$ let
\begin{align*}
\xi^\pm_{t,\ell,\vec\eta}(\cdot)
&\equiv\nu_{G^k_{2n}}(\si_{i^\ell_\pm}=\cdot\mid\vec\si_{S_t}=\vec\eta),\\
\ze^\pm_{t,\ell,\vec\eta}(\cdot)
&\equiv\nu_{G^k_{2n}}(\si_{j^\ell_\pm}=\cdot\mid\vec\si_{S_t}=\vec\eta),
\end{align*}
so that
$$
\nu_{G_{2n}}[(\si_{i^\ell_+},\si_{j^\ell_-})=\cdot
	\mid \vec\si_{S_t}=\vec\eta]
= \xi^+_{t,\ell,\vec\eta} \otimes_\psi \ze^-_{t,\ell,\vec\eta}.$$
By Thm.~\ref{t:loc}, the conditional measures $\nu^+_{G_{2n}}(\vec\si_{\ball{t}{i^\ell_+}}=\cdot)$ converge to $\mu^+$. But by \eqref{e:magnetization}, $Y(\vec\si)$ agrees with $Y_t(\vec\si)\equiv \sgn \sum_{i\in V\setminus B_t}\col_i\si_i$ with high probability, so that convergence also holds if we replace $\nu^+_{G_{2n}}$ by $\nu^{\pm t}_{G_{2n}}(\cdot)\equiv\nu_{G_{2n}}(\cdot\mid Y_t(\vec\si)=\pm)$. In particular,
\begin{align*}
0&=\lim_{t\to\infty}\lim_{n\to\infty}
\E_{2n}\Big[
\tvdb{
\sum_{\vec\eta}\nu^{+t}_{G_{2n}}(\vec\si_{S_t}=\vec\eta)
\xi^+_{t,\ell,\vec\eta}\otimes_\psi \ze^-_{t,\ell,\vec\eta}
- h^+ \otimes_\psi h^-
}\Big]\\
&=\lim_{t\to\infty}\lim_{n\to\infty}
\E_{2n}\Big[
\tvdb{
\angl{\xi^+_{t,\ell,\vec\si_{S_t}}
\otimes_\psi \ze^-_{t,\ell,\vec\si_{S_t}}}_{\nu^{+t}_{G_{2n}}}
- h^+ \otimes_\psi h^-}\Big]
\end{align*}
On the other hand, it is easily seen that $(h\otimes_\psi h')(1,0)$ is maximized by taking $h(1)$ and $h'(0)$ as large as possible. But in the limit $t\to\infty$ the values $\xi^\pm_{t,\ell,\vec\eta}(1),\ze^\pm_{t,\ell,\vec\eta}(1)$ (with $\vec\eta$ arbitrary) are sandwiched between $h^\pm(1)$, so it must be that
\beq\label{e:xiLimit}
0=\lim_{t\to\infty}\lim_{n\to\infty}
\E_{2n}\Big[ \anglb{
\tvd{\xi^+_{t,\ell,\vec\si_{S_t}}-h^+}
+\tvd{\ze^-_{t,\ell,\vec\si_{S_t}}-h^-}
}_{\nu^{+t}_{G_{2n}}} \Big].
\eeq

We now claim that \eqref{e:xiLimit} continues to hold after removal of the edges $(i^\ell_\pm,j^\ell_\mp)$. Indeed,
\beq\label{e:ZComparison}
\f{\nu^{+t}_{G_{2n}}(\vec\si_{S_t}=\vec\eta)}{\nu^{+t}_{G^k_{2n}}(\vec\si_{S_t}=\vec\eta)}
= \f{Z^{+t}_\OUT(\vec\eta) Z_\IN(\vec\eta)}
{Z^{+t}_\OUT(\vec\eta) Z^k_\IN(\vec\eta)}
\cdot \f{\sum_{\vec\eta'} Z^{+t}_\OUT(\vec\eta') Z^k_\IN(\vec\eta')}
{\sum_{\vec\eta'} Z^{+t}_\OUT(\vec\eta') Z_\IN(\vec\eta')}
\eeq
where
\begin{align*}
Z^{\pm t}_\OUT(\vec\eta)
&\equiv Z_{G_{2n} \setminus B_{t-1}}
[\set{
\vec\si_{G_{2n}\setminus B_{t-1}}
: Y_t(\vec\si)=\pm \text{ and } \vec\si_{S_t}=\vec\eta}],\\
Z_\IN(\vec\eta)
&\equiv Z_{B_t}
	[\set{\vec\si_{B_t} : \vec\si_{S_t}=\vec\eta}],\\
Z^k_\IN(\vec\eta)
&\equiv Z_{B_t \cap G^k_{2n}}
	[\set{\vec\si_{B_t} : \vec\si_{S_t}=\vec\eta}].
\end{align*}
Now note that for $k$ bounded and $t$ large we have $Z_\IN(\vec\eta)\asymp Z^k_\IN(\vec\eta)$ uniformly over $\vec\eta$: for Ising interactions at non-zero temperature this is obvious, while for the hard-core model
$$\f{Z_\IN(\vec\eta)}{Z^k_\IN(\vec\eta)}
= \prod_{\ell=1}^k
\ls
[1-\xi^+_{t,\ell,\vec\eta}(1)\ze^-_{t,\ell,\vec\eta}(1)]
[1-\xi^-_{t,\ell,\vec\eta}(1)\ze^+_{t,\ell,\vec\eta}(1)]
\rs$$
which for $t$ large is $\asymp 1$ uniformly over $\vec\eta$. Since the $\xi^\pm_{t,\ell,\vec\eta}$ and $\ze^\pm_{t,\ell,\vec\eta}$ are $\vec\eta$-measurable, it follows from \eqref{e:ZComparison} that \eqref{e:xiLimit} continues to hold with $\nu^{+t}_{G^k_{2n}}$ in place of $\nu^{+t}_{G_{2n}}$. Since the spins $(\si_w)_{w\in W}$ are independent under $\nu^{\pm t}_{G^k_{2n}}(\cdot\mid\vec\si_{S_t})$, this further implies
\beq\label{e:conv.to.product}
0
=\lim_{t\to\infty}\lim_{n\to\infty}
\E_{2n}\Big[
\tvd{\nu^{+t}_{G^k_{2n}}(\vec\si_W=\cdot)
-Q^+_W}\Big].
\eeq
Finally, by a similar argument as before $\lim_{n\to\infty}\nu_{G^k_{2n}}(Y(\vec\si)=Y_t(\vec\si))=1$, so \eqref{e:conv.to.product} holds with $\nu^+_{G^k_{2n}}$ in place of $\nu^{+t}_{G^k_{2n}}$ which gives the result.
\epf
\eppn

We now demonstrate how to use Propn.~\ref{p:product} to establish a randomized reduction from approximating the partition function to the problem of approximate $\maxcut$ on 3-regular graphs, which is $\NP$-hard~\cite{\AKMaxcut}.

Let $H$ be a 3-regular graph on $m$ vertices and construct the bipartite graph $G=G^{3 k}_{2n}$ by the procedure described above. By Lem.~\ref{l:zero} and Propn.~\ref{p:product}, for any $\ep>0$ there exists $n(\ep)$ large enough such that the following hold with positive probability:
\bnm[{\sc (i)}]
\item $G^{3k}_{2n}$ was formed by removing $3k$ distinct edges from a $d$-regular graph $G_{2n}$;
\item $\nu_{G^{3k}_{2n}}(Y(\vec\si)=+)\le (1+\ep)/2$; and
\item $\nu^\pm_{G^{3k}_{2n}}(\vec\si_W) / Q^\pm_W(\vec\si_W)
\in[1-\ep,1+\ep]$ for all $\vec\si_W$.
\enm
Consequently, for given $\ep$ we may find $G^{3k}_{2n}$ satisfying properties {\sc (i)-(iii)} within finite time by deterministic search. We then construct from $H$ and $G$ a new graph $H^G$ as follows:
\bitm
\item For each vertex $x\in H$ let $G_x$ be a copy of $G$, and denote by $W^\pm_x$ the vertices of $G_x$ corresponding to $W^\pm$ in $G$. Let $\wh H^G$ be the disjoint union of the $G_x$, $x\in H$.
\item For every edge $(x,y)\in H$, add $2k$ edges between $W^+_x$ and $W^+_y$ and similarly $2k$ edges between $W^-_x$ and $W^-_y$. This can be done deterministically in such a way that the resulting graph, which we denote $H^G$, is $d$-regular.
\eitm
We write a spin configuration on $\wh H^G$ or $H^G$ as $\vec\si\equiv(\vec\si_x)_{x\in H}$ where $\vec\si_x$ is the restriction of $\vec\si$ to $G_x$. We write $Y_x\equiv Y(\vec\si_x)$ for the phase of each $\vec\si_x$, and $\cY(\vec\si)\equiv(Y(\vec\si_x))_{x\in H}\in\set{0,1}^H$. Write $Z_{H^G}(\cY)$ for the partition function for the two-spin model on $H^G$ restricted to configurations of phase $\cY$, and define likewise $Z_{\wh H^G}(\cY)$.

Recalling \eqref{e:z.psi.h}, let
$$\Gam\equiv
z(h^+\otimes_\psi h^+)
z(h^-\otimes_\psi h^-),\quad
\Thet\equiv z(h^+\otimes_\psi h^-)^2,$$
and note that for anti-ferromagnetic two-spin models in non-uniqueness regimes, $\Thet>\Gam$.

\blem\label{l:HGcutProb}
For $G$ satisfying properties {\sc (i)-(iii)},
$$[(1-\ep)/2]^m
\le\f{ Z_{H^G}/Z_{\wh H^G} }
{\Gam^{2k\abs{E(H)}} (\Thet/\Gam)^{2k\,\maxcut(H)}}
\le(1+\ep)^m.$$

\bpf
By~{\sc (ii)},
\beq\label{e:ratio.zy.z}
(1-\ep)^m
\le 2^m \f{Z_{\wh H^G(\cY)}}{Z_{\wh H^G}}
\le (1+\ep)^m
\eeq
for all $\cY\in\set{0,1}^H$. By {\sc (iii)}, the ratio
$$\f{Z_{H^G}(\cY)}{Z_{\wh H^G}(\cY)}
=\sum_{x\in H} \sum_{\vec\si_{W_x}} \nu^{Y_x}_{G_x}(\vec\si_{W_x})
\prod_{(i,j)\in E(H^G)\setminus E(\wh H^G)}
\psi(\si_i,\si_j)$$
is within a $(1\pm\ep)^m$ factor of
$$\sum_{x\in H} \sum_{\vec\si_{W_x}}
Q^{Y_x}(\vec\si_{W^+_x})
\prod_{(i,j)\in E(H^G)\setminus E(\wh H^G)}
\psi(\si_i,\si_j),$$
which by direct calculation equals
$$\Gam^{2k\abs{E(H)}} (\Thet/\Gam)^{2k\,\cut(\cY)}$$
where $\cut(\cY)\equiv\abs{\{(x,y)\in E(H): Y_x\neq Y_y\}}$, the number of edges crossing the cut of $H$ induced by $\cY$. Combining with \eqref{e:ratio.zy.z} gives
$$Z_{H^G}
= \sum_\cY \f{Z_{H^G}(\cY)}{Z_{\wh H^G}(\cY)}
	Z_{\wh H^G}(\cY)
\le (1+\ep)^{2m}
	\Gam^{2k\abs{E(H)}}(\Thet/\Gam)^{2k\,\maxcut(H)}
	Z_{\wh H^G}$$
and similarly
$$Z_{H^G}
\ge 2^{-m}(1-\ep)^{2m} \Gam^{2k\abs{E(H)}}(\Thet/\Gam)^{2k\,\maxcut(H)} Z_{\wh H^G}.$$
Rearranging gives the stated result.
\epf
\elem

Using this lemma we now complete the reduction to approximate $\maxcut$:

\begin{proof}[Proof of Thms.~\ref{t:hc.comp} and~\ref{t:ising.comp}]
Let $H$ be a $3$-regular graph on $m$ vertices, and note that the maximum cut of $H$ is at least $3m/4$, the expected value of a random cut. Construct $\wh H^G$, $H^G$ as above. Since $\wh H^G$ is a disjoint collection of constant-size graphs, its partition function can be computed in polynomial time. Suppose $Z_{H^G}$ could be approximated within a factor of $e^{c\abs{H^G}}$ in polynomial time for any $c>0$: rearranging the result of Lem.~\ref{l:HGcutProb} gives
\beq\label{e:maxcut.bounds}
\f{\displaystyle \log \lp
	\f{Z_{H^G}/Z_{\wh H^G}}
	{\Gam^{2k\abs{E(H)}} (1+\ep)^m}\rp}{2k\log(\Thet/\Gam)}
\le\maxcut(H)
\le \f{\displaystyle \log\lp
	\f{Z_{H^G}/Z_{\wh H^G}}
	{\Gam^{2k\abs{E(H)}} [(1-\ep)/2]^m}\rp}{2k\log(\Thet/\Gam)},
\eeq
so within polynomial time one obtains upper and lower bounds for $\maxcut(H)$ which differ by $O[(c\abs{G}+1)m/k]$. Taking $k$ large and $c$ small then allows to compute $\maxcut(H)$ up to an arbitrarily small multiplicative error: that is, we have completed the reduction to a $\pras$ for $\maxcut$ on $3$-regular graphs, in contradiction of the result of~\cite{\AKMaxcut}.
\epf

\subsection*{Acknowledgements}

We thank Andreas Galanis, Daniel \v{S}tefankovi\v{c}, and Eric Vigoda for describing to us their methods and for sending us a draft of their paper. We thank
Amir Dembo, David Gamarnik, Andrea Montanari, Alistair Sinclair, Piyush Srivastava, and David Wilson for helpful conversations.

\bibliography{refs}
\bibliographystyle{abbrv}

\end{document}